\documentclass[11pt]{amsart}

\usepackage[title]{appendix}
\usepackage{titlesec}
\usepackage{amsmath,amssymb,mathtools, amsfonts}
\usepackage{mathrsfs}
\usepackage{graphicx}
\usepackage{enumitem}
\usepackage{microtype}
\usepackage[LGR,T1]{fontenc}
\usepackage{alphabeta}
\usepackage{amscd}
\usepackage[colorlinks=true,
            linkcolor=blue,
            citecolor=blue,
            urlcolor=blue]{hyperref}
\usepackage[nameinlink,noabbrev]{cleveref}
\usepackage[margin=1in]{geometry}

\AtBeginDocument{%
}

\titleformat{\section}
  {\normalfont\Large\bfseries}
  {\thesection}{1em}{}
  [\titlerule]

\newtheorem{theorem}{Theorem}[section]

\newtheorem{lemma}[theorem]{Lemma}
\newtheorem{corollary}[theorem]{Corollary}

\theoremstyle{definition}

\theoremstyle{remark}

\theoremstyle{plain}

\newtheorem*{repthm}{\repthmname}
\newcommand{\repthmname}{}
\newenvironment{restate}[2][Theorem]
  {\renewcommand{\repthmname}{#1~\ref*{#2}}\begin{repthm}}
  {\end{repthm}}

\def\inmod#1{\allowbreak\mkern5mu({\operatorname{mod}\,#1)}}

\title{On p-Group Actions and Index Divisibilities}

\author{Kyle Gettig}
\address{University of South Florida}
\email{kgettig@usf.edu}


\begin{document}

\begin{abstract}
We prove a mod-$ p $ divisibility result for indices of elliptic operators assuming a smooth action by a finite $ p $-group on the base manifold of dimension $ n\le 2(p-1) $ with no stationary points and which lifts equivariantly to the operator bundles.

We show as one consequence that if such manifolds are spin or stably almost complex, with the group action lifting to the relevant structure, then their class in the relevant bordism theory modulo $ p $ is fully determined by one genus, and indeed vanishes for $ n<2(p-1) $.

We show as another consequence that there exists no closed connected $ 24 $-dimensional spin manifold $ M $ satisfying $ {\hat{A}(M)=1}, {\hat{A}(M;T_\mathbb{C}M)=0} $, and for which there exists a smooth non-trivial action by the monster group, resolving a question of Hirzebruch.
\end{abstract}

\maketitle

\section{Introduction and Statement of Results}
Given a group $ G $, the question of which manifolds are acted on non-trivially by $ G $ is a central structural one in topology. Assumptions on the particular group $ G $ (e.g. abelian, finite of some order, simple or semi-simple, etc.), additional structure on the manifold (e.g. some $ H $-structure), and/or specific properties of the action (e.g. smooth, effective, Hamiltonian, preserves the aforementioned manifold structure, etc.) are frequently combined to show obstructions to the existence of such actions, often in terms of numerical invariants. As a classical example, Poincaré-Hopf implies that the Euler characteristic is an obstruction to a free action by $ S^1 $ on a closed manifold. Similar results proliferated with the development of equivariant index theories such as Atiyah-Singer in the latter half of the 20th century \cite{borel1960seminar, atiyah1968index, bredon1972introduction, schultz1985group}.

Our main theorem is a result in this category:
\begin{restate}{thm:main}
Let $ p $ be an odd prime and $ P $ be a finite $ p $-group, and let $ M $ be a closed manifold of dimension $ n $ on which $ P $ acts smoothly without stationary points. Let $ D $ be a $ P $-equivariant elliptic operator over $ M $.
\begin{enumerate}[label=(\alph*)]
\item If $ n<2(p-1) $, then $ \mathrm{ind}(D)\equiv 0\inmod{p} $.
\item If $ n=2(p-1) $ and $ M $ is connected, then $ \mathrm{ind}(D\otimes W)\equiv \mathrm{rk}(W)\,\mathrm{ind}(D)\inmod{p} $ for every $ P $-equivariant vector bundle $ W $ on which the kernel of $ P\to\mathrm{Diff}(M) $ acts trivially.
\end{enumerate}
\end{restate}
As a consequence, for spin and complex bordism whose classes are known to be determined by characteristic numbers, we have the following bordism-level corollary.
\begin{restate}[Corollary]{thm:divisibility-of-characteristic-numbers}
Let $ P $ be a $ p $-group for odd $ p $ which acts smoothly on a closed spin (resp. stably almost complex) manifold $ M $ of real dimension $ n=2m $ without stationary points. In the complex case, assume $ P $ preserves the stably almost complex structure.
\begin{enumerate}[label=(\alph*)]
\item If $ m<p-1 $, then the bordism class $ [M] $ is divisible by $ p $ in spin (resp. complex) bordism.
\item If $ m=p-1 $ in the complex case, then $ [M]\equiv \mathrm{Td}(M)[\mathbb{C}P^{p-1}]\inmod{p} $ in complex bordism.
\item If $ m=p-1 $ in the spin case, then $ [M]\equiv 2^{-1}\hat{A}(M)[V_p]\inmod{p} $ in spin bordism, where $ V_p $ is a smooth hypersurface of degree $ 2p $ in $ \mathbb{C}P^p $.
\end{enumerate}
\end{restate}
These together form a smooth generalization of cases (iii) of the results of Haution \cite{haution2019fixed}, which consider only projective varieties and their algebraic invariants. Haution's results build on earlier work of tom Dieck \cite{dieck1970actions}, who produced the corresponding results in complex bordism for abelian $ P $.

As a near-immediate consequence, we get the following result for spin manifolds.

\begin{restate}[Corollary]{thm:main-thm-for-spin}
Let $ P $ be a $ p $-group for odd $ p $ which acts smoothly and without stationary points on a closed connected spin manifold $ M $ of dimension $ 2m $.
\begin{enumerate}[label=(\alph*)]
\item If $ m<p-1 $, then $ \hat{A}(M;W)\equiv 0\inmod{p} $ for every $ P $-equivariant complex vector bundle $ W $.
\item If $ m=p-1 $, then $ \hat{A}(M;W)\equiv\mathrm{rk}(W)\hat{A}(M)\inmod{p} $ for every $ P $-equivariant complex vector bundle $ W $ on which the kernel of $ P\to\mathrm{Diff}(M) $ acts trivially. In particular, $ \mathrm{sign}(M)\equiv\hat{A}(M)\inmod{p} $.
\end{enumerate}
\end{restate}

As a particular application, this implies a negative answer to a question of Hirzebruch.

\begin{restate}[Corollary]{thm:moonshine}
No closed, connected $ 24 $-dimensional spin manifold $ M $ with $ \hat{A}(M)=1 $ and $ \hat{A}(M;T_\mathbb{C}M)=0 $ can be acted on smoothly and non-trivially by the monster group.
\end{restate}

The genesis of this question lies in so-called `monstrous moonshine': the observation in the late 1970s of several numerical `coincidences' between the theory of modular forms, and in particular the j function, and representations of the (then still conjectured) monster group \cite{conway1979monstrous}. Much of the connection between these two was clarified by the discovery of the `moonshine module' in the late 80s by Frenkel, Lepowsky, and Meurman: a vertex operator algebra on which the monster acts as the automorphism group \cite{frenkel1989vertex}.

Soon after, Hirzebruch posed his `prize question': does there exist a $ 24 $-dimensional spin manifold $ M $ with $ p_1(M)=0 $, $ \hat{A}(M)=1 $, and $ \hat{A}(M;T_\mathbb{C}M)=0 $ on which the monster acts by diffeomorphisms \cite[p.87]{hirzebruch1992manifolds}? Firstly, the Witten genus of such a manifold would be $ (j-744)\Delta $, where $ \Delta $ is the modular discriminant; this would then provide a direct link between $ j $ and representations of the monster via equivariant twisted Dirac indices. Secondly, from a physical perspective, a non-linear sigma model with target $ M $ could conceivably produce a conformal field theory whose chiral algebra is exactly the moonshine module, giving an appealing differential-geometric motivation for the origin of the latter (see, e.g. \cite{ganter2007hecke}).

The existence of such a manifold without the monster action was demonstrated by Hopkins and Mahowald in 2002 using techniques of homotopy theory \cite{mahowald2002structure}. Our corollary shows that no such manifold can have a non-trivial smooth monster action, resolving Hirzebruch's question in the negative.

\vspace{10pt}

The paper is organized as follows: We first state and prove two lemmas necessary for our main theorem. We then state and prove our main theorem. We then state and prove several corollaries to our main theorem, as presented in the introduction.

\vspace{10pt}

Claude Opus 5.5 was used in the development of this paper. All writing is the author's own, and the author has exercised all due diligence in its use.

\section{Two Lemmas}
In the sequel, let $ p $ be an odd prime, $ C_p $ be the cyclic group of order $ p $, and $ \zeta=e^{2\pi i/p} $.

\begin{lemma}\label{thm:index-lefschetz-equivalence}
Let $ D: E\to F $ be an elliptic operator over the manifold $ M $. Suppose $ C_p $ acts smoothly on $ M $, and this action lifts to $ E $ and $ F $ such that it commutes with $ D $. Then the Lefschetz number $ L_g(D) $ for any non-identity $ g\in C_p $ lives in $ \mathbb{Z}[\zeta] $, and
\[ \mathrm{ind}(D)\equiv L_g(D)\inmod{(1-\zeta)}. \]
\end{lemma}

\begin{proof}
This follows immediately from all eigenvalues of $ g $ being of the form $ \zeta^k\equiv 1\inmod{(1-\zeta)} $, whence traces become dimensions when taking mod $ (1-\zeta) $.
\end{proof}

Let $ g\in C_p $ act on an $ n $-dimensional closed manifold $ M $, and let $ D $ be a $ C_p $-equivariant elliptic operator over $ M $ with symbol $ \sigma(D)\in K_{C_p}(TM) $. Let $ F $ be a component of the fixed point set $ M^g $ of dimension $ 2d_F $ (we permit $ d_F $ to be half-integer), and let $ F $ have normal bundle $ N $ of complex dimension $ r_F $, so that $ d_F+r_F=n/2 $. The action of $ g $ on $ N $ lets us split $ N\cong\oplus_mN_m $, where each $ N_m $ is a complex vector bundle on which $ g $ acts as multiplication by $ \zeta^m $.

One form of the fixed-point formula of Atiyah-Segal \cite{atiyah1968indexII} is that $ L_g(D)=\sum_F c_F $, where
\begin{equation}\label{eq:atiyah-segal-term}
c_F=\text{ind}\left(\frac{i_F^*\sigma(D)(g)}{\lambda_{-1}(N_\mathbb{C})(g)}\right).
\end{equation}
By the Thom isomorphism of equivariant K-theory (see \cite{atiyah1968bott}), $ i^*_F\sigma(D)=\tau_F\cdot\lambda_{-1}(N) $ for a unique $ \tau_F\in K_{C_p}(TF)$. Since $ N_\mathbb{C}\cong N\oplus\overline{N} $, a cancellation of $ \lambda_{-1}(N) $ in (\ref{eq:atiyah-segal-term}) gives
\begin{equation}\label{eq:atiyah-segal-term-2}
c_F=\text{ind}(\tau_F(g)\cdot E_F^{-1}),
\end{equation}
where $ E_F=\lambda_{-1}(N^*)(g) $ is the equivariant K-theory Euler class.

We say that $ V\in K(F)\otimes\mathbb{Q} $ has degree at least $ j $ if the Chern character $ \text{ch}(V) $ lives in cohomological degrees at least $ 2j $. As $ g $ acts trivially on $ TF $, $ \tau_F $ splits by eigenvalues as $ \tau_F(g)=\sum_k\zeta^k\tau_{F,k} $. Each $ \tau_{F,k} $ lies in $ K(TF) $, and so is the symbol of some elliptic operator $ D_{F, k} $, which is an operator induced by $ D $ over $ F $.

\begin{lemma}\label{thm:form-of-c-F}
For each component $ F $ of $ M^g $ there are classes $ V_0, \cdots, V_{\lfloor d_F\rfloor} $ in $ K(F)\otimes\mathbb{Z}[\zeta] $, built from the $ N_m $ by exterior powers independently of $ D $ and with $ V_j $ of degree at least $ j $, such that
\[ c_F=\displaystyle\sum_k\displaystyle\sum_{j=0}^{\lfloor d_F\rfloor}\zeta^k(1-\zeta)^{-(r_F+j)}\mathrm{ind}(D_{F,k}\otimes V_j). \]
\end{lemma}

\begin{proof}
This follows relatively directly from the formula (\ref{eq:atiyah-segal-term-2}) for $ c_F $ and the expansion
\[ E_F=\lambda_{-1}(N^*)(g)=\displaystyle\prod_m\displaystyle\sum_k(-\zeta^{-m})^k\Lambda^kN^*_m. \]
Fix an $ m $ and let $ r_m=\text{rk}_\mathbb{C}N_m $. A direct computation gives
\[ \lambda_{-1}(N_m^*)(g)=(1-\zeta^{-m})^{r_m}(1+X_m), \]
where
\begin{equation}\label{eq:Xm-expansion}
X_m=\displaystyle\sum_{j\ge 1}(1-\zeta^{m})^{-j}\gamma^j(N_m^*-r_m).
\end{equation}
Here $ \gamma^j $ is the $ j $th $ \gamma $-operation; we may split $ N_m^*-r_m $ as a sum $ \oplus_i(L_i-1) $ for complex line bundles $ L_i $, and in this case $ \gamma^j(N_m^*-r_m) $ is the $ j $th elementary symmetric function of the $ L_i-1 $. In particular, $ \gamma^j(N_m^*-r_m) $ has degree at least $ j $.

As $ X_m $ has rank $ 0 $, it's nilpotent, so $ 1+X_m $ is invertible in K-theory; we can thus express
\begin{equation}\label{eq:E_F-inverse-expansion}
\begin{aligned}
E_F^{-1}&=\displaystyle\prod_m\lambda_{-1}(N_m^*)(g)^{-1} \\ &=\displaystyle\prod_m(1-\zeta^{-m})^{-r_m}(1+X_m)^{-1}.
\end{aligned}
\end{equation}
We expand this product together with (\ref{eq:Xm-expansion}) and, for each fixed $ j $, collect together terms featuring $ \gamma^{a_1}\cdots\gamma^{a_l} $ for $ a_1+\cdots+a_l=j $, pull out all possible factors of $ (1-\zeta)^{-1} $, and let the result be our $ V_j $.

Finally, we count how many such factors of $ (1-\zeta)^{-1} $ get so attached to each $ V_j $. The complex coefficients in each term of (\ref{eq:E_F-inverse-expansion}) lie in $ (1-\zeta)^{-r_m}\mathbb{Z}[\zeta] $, so we get $ (1-\zeta)^{-r_F} $ for free independently of $ j $. We further get $ (1-\zeta)^{-j} $ attached to each $ \gamma^j $ term; this together gives the $ (1-\zeta)^{-(r_F+j)} $ factor in the lemma.
\end{proof}

\section{Main Theorem}

\begin{theorem}\label{thm:main}
Let $ p $ be an odd prime and $ P $ be a finite $ p $-group, and let $ M $ be a closed manifold of dimension $ n $ on which $ P $ acts smoothly without stationary points. Let $ D $ be a $ P $-equivariant elliptic operator over $ M $.
\begin{enumerate}[label=(\alph*)]
\item If $ n<2(p-1) $, then $ \mathrm{ind}(D)\equiv 0\inmod{p} $.
\item If $ n=2(p-1) $ and $ M $ is connected, then $ \mathrm{ind}(D\otimes W)\equiv \mathrm{rk}(W)\,\mathrm{ind}(D)\inmod{p} $ for every $ P $-equivariant vector bundle $ W $ on which the kernel of $ P\to\mathrm{Diff}(M) $ acts trivially.
\end{enumerate}
\end{theorem}
We make two quick remarks on the $ W $ in part (b). Firstly, the condition on $ W $ is immediately satisfied for any $ W $ built naturally from the tangent bundle $ TM $ (i.e. tensor, symmetric, exterior powers, and their sums and differences). Secondly, once part (b) is true, it immediately extends to $ P $-equivariant virtual K-theory bundles as well.

\begin{proof}
Let $ O=\mathbb{Z}_{(p)}[\zeta] $, the ring $ \mathbb{Z}[\zeta] $ with every prime other than $ p $ inverted. It has unique maximal ideal $ \mathfrak{p}=(1-\zeta) $ and a $ p $-adic valuation $ v $ which we renormalize to $ v(1-\zeta)=1 $, so that $ v(p)=p-1 $. Standard algebra says that
\[ \mathfrak{p}=\{x\in O\,|\,v(x)>0\}. \]
The key fact we abuse is that $ \mathfrak{p}\cap\mathbb{Z}=p\mathbb{Z} $, so any integers in $ O $ with positive valuation are necessarily divisible by $ p $.

\emph{(a)}\;We induct on pairs $ (|P|, n) $, ordered lexicographically. $ P $ trivial is vacuous by the stationary-point-free assumption. We may assume $ M $ is connected, as components of $ M $ in the same orbit of $ P $ contribute equally to the operator index by equivariance, and such orbits are of size either $ 1 $ or some multiple of $ p $.

Let $ z\in P $ be a central element of order $ p $. We split into cases where $ z $ acts non-trivially on $ M $ and where $ z $ acts trivially.

Suppose firstly that $ z $ acts non-trivially. By Lemma \ref{thm:index-lefschetz-equivalence}, $ \text{ind}(D)\equiv L_z(D)=\sum_Fc_F\inmod{\mathfrak{p}} $, summing over fixed components $ F $ of $ M^z $. Since $ z $ is central, the full group $ P $ permutes the components $ F $, and by naturality $ c_{gF}=c_F $ for $ g\in P $. By Lemma \ref{thm:form-of-c-F}, $ v(c_F)\ge-\lfloor n/2\rfloor $.

\begin{itemize}
\item If $ F $ is moved by $ P $, the orbit size is a multiple of $ p $, and so the orbit contributes to $ L_z(D) $ terms with total valuation no less than $ p-1-\lfloor n/2\rfloor>0 $.
\item If $ F $ is invariant under $ P $, induction on the dimension shows that all indices in the summed expression in Lemma \ref{thm:form-of-c-F} are divisible by $ p $ (full $ P $-equivariance of the restricted operators, which is necessary for the induction, follows from centrality of $ z $). The terms summing to $ c_F $ therefore each have valuation no less than $ p-1-\lfloor n/2\rfloor >0 $, and so the same is true of $ c_F $ itself.
\end{itemize}
Therefore $ L_z(D) $ is a sum of terms of positive valuation, implying $ L_z(D)\in\mathfrak{p} $, and hence $ \text{ind}(D)\in\mathfrak{p} $ by Lemma \ref{thm:index-lefschetz-equivalence}. We conclude that $ \text{ind}(D)\equiv 0\inmod{p} $.

Suppose next that $ z $ acts trivially on $ M $. By the hypothesis that $ D $ is $ P $-equivariant, the action of $ z $ splits by eigenvalues $ D=\oplus_k D_k $ with each $ D_k $ elliptic, so that $ z $ acts on the equivariant index $ \text{ind}_P(D_k) $ as the scalar $ \zeta^k $. We split into three further cases.

\begin{itemize}
\item If $ k=0 $ so that $ z $ acts fully trivially, all actions descend to actions of $ P/\langle z\rangle $ which preserves equivariance. Since $ z $ already fixed everything, the $ P/\langle z\rangle $ action has no stationary points. By induction on the size of the $ p $-group, $ \text{ind}(D_0)\equiv 0\inmod{p} $.
\item If $ k\not=0 $ and $ z $ is in the commutator subgroup $ [P, P] $, then all one-dimensional representations are acted on trivially by $ z $, meaning no one-dimensional representations can appear in $ \text{ind}_P(D_k) $ by assumption that $ z $ acts as $ \zeta^k $. All other irreducible representations have dimensions that are positive powers of $ p $, and so $ \text{ind}(D_k)\equiv 0\inmod{p} $.
\item If $ k\not=0 $ and $ z $ is not in the commutator subgroup $ [P, P] $, then there exists some one-dimensional character $ \chi $ of $ P $ with $ \chi(z)=\zeta^{-k} $. The twisted operator $ D_k\otimes\chi $ is $ P $-equivariant and has the same index as $ D_k $ but is acted on by $ z $ with weight $ 0 $, and so the first case applies, giving $ \text{ind}(D_k)\equiv 0\inmod{p} $.
\end{itemize}
Summing over $ k $ therefore gives $ \text{ind}(D)\equiv 0\inmod{p} $, as desired. \\

\emph{(b)}\;We induct on $ |P| $. The shape of this proof will largely mirror the shape of the proof of part (a); as in that part, we fix an element $ z\in P $ of order $ p $ in the center. Let $ R=\text{rk}(W) $. Unlike in part (a), the kernel condition on the action of $ P $ on $ W $ prevents a reduction to the connected case, so we hypothesize a priori $ M $ is connected.

The case in which $ z $ acts trivially on $ M $ is a near-direct echo of the corresponding case in part (a), so we begin with it, being somewhat briefer. By assumption $ z $ also acts trivially on $ W $, and by equivariance, we may again split $ D=\oplus_k D_k $ by eigenvalues of the action of $ z $.
\begin{itemize}
\item If $ k=0 $, $ P/\langle z\rangle $ acts on all relevant objects without stationary points, and the kernel of $ P/\langle z\rangle\to\text{Diff}(M) $ still acts trivially. Induction on the group size therefore gives $ \text{ind}(D_0\otimes W)\equiv R\,\text{ind}(D_0)\inmod{p} $.
\item If $ k\not=0 $ and $ z\in [P, P] $, then $ \text{ind}(D_k\otimes W)\equiv\text{ind}(D_k)\equiv 0\inmod{p} $ completely analogously to the part (a) case.
\item If $ k\not=0 $ and $ z\not\in [P, P] $, then we can again twist $ D_k $ by a character $ \chi $ with $ \chi(z)=\zeta^{-k} $ and return to the $ k=0 $ case, which is proven.
\end{itemize}
Summing over $ k $ therefore gives the desired $ \text{ind}(D\otimes W)\equiv R\,\text{ind}(D)\inmod{p} $.

Now suppose $ z $ acts non-trivially on $ M $. The failure of our ability to repeat the corresponding argument in part (a) verbatim comes precisely at the $ j=d_F $ term in the summation formula for $ c_F $ in Lemma \ref{thm:form-of-c-F}; this is the point at which we have one too many factors of $ (1-\zeta)^{-1} $, and we can no longer guarantee a positive valuation. If we let $ c_F(W) $ denote the term $ c_F $ associated to the fixed point component $ F $ in $ M^z $ and the operator $ D\otimes W $, and similarly let $ c_F(1) $ denote the same corresponding to the operator $ D $, what we need is an extra factor of $ 1-\zeta $ to appear in the $ j=d_F $ terms of the Lemma \ref{thm:form-of-c-F} summation expansion of
\[ c_F(W)-Rc_F(1). \]

This is where we need the fact that $ V_{d_F} $ has degree at least $ d_F $. The relevant part of $ c_F(W)-Rc_F(1) $ that we're looking at is
\[ \displaystyle\sum_k\zeta^k\left(\displaystyle\sum_m\zeta^m\text{ind}(D_{F,k}\otimes (W|_F)_m\otimes V_{d_F})-R\,\text{ind}(D_{F,k}\otimes V_{d_F})\right), \]
where $ W|_F=\oplus_m(W|_F)_m $ splits via eigenvalues of $ z $. As the cohomological degree of $ \text{ch}(V_{d_F}) $ saturates the dimension of $ F $, by Atiyah-Singer,
\[ \text{ind}(D_{F,k}\otimes (W|_F)_m\otimes V_{d_F})=\text{rk}((W|_F)_m)\text{ind}(D_{F,k}\otimes V_{d_F}). \]
The above sum therefore becomes
\[ \left(\displaystyle\sum_k \zeta^k\text{ind}(D_{F,k}\otimes V_{d_F})\right)\left(\displaystyle\sum_m\zeta^m\text{rk}((W|_F)_m)-R\right). \]
The second factor gives us the extra factor of $ 1-\zeta $ that we need at $ j=d_F $. With this extra factor, the valuation argument of part (a) goes through.
\end{proof}

\section{Consequences}
\begin{corollary}\label{thm:divisibility-of-characteristic-numbers}
Let $ P $ be a $ p $-group for odd $ p $ which acts smoothly on a closed spin (resp. stably almost complex) manifold $ M $ of real dimension $ 2m $ without stationary points. In the spin case, assume $ M $ is connected; in the complex case, assume $ P $ preserves the stably almost complex structure.
\begin{enumerate}[label=(\alph*)]
\item If $ m<p-1 $, then the bordism class $ [M] $ is divisible by $ p $ in spin (resp. complex) bordism.
\item If $ m=p-1 $ in the complex case, then $ [M]\equiv \mathrm{Td}(M)[\mathbb{C}P^{p-1}]\inmod{p} $ in complex bordism.
\item If $ m=p-1 $ in the spin case, then $ [M]\equiv 2^{-1}\hat{A}(M)[V_p]\inmod{p} $ in spin bordism, where $ V_p $ is a smooth hypersurface of degree $ 2p $ in $ \mathbb{C}P^p $.
\end{enumerate}
\end{corollary}
\begin{proof}
Remark that $ V_p $ is simply connected in (c) and therefore has unique spin structure. We assume, by considering connected components, that $ M $ is connected. We can give a unified proof of all three parts. Firstly, suppose all the KO-theory (resp. K-theory) characteristic numbers of a bordism class vanish mod $ p $. Then the rational class $ \frac{1}{p}[M] $ in the relevant rational bordism ring has its characteristic numbers also satisfying Stong's characteristic number relations \cite{stong1965relations, stong1966relations} (see also \cite{hattori1966integral}). In the complex case, this immediately associates $ \frac{1}{p}[M] $ with some class in $ \Omega^U_{2m} $. In the spin case, this only necessarily associates it with some class in $ \Omega^\text{Spin}_{2m}/\text{Torsion} $; we then need the additional fact that torsion in spin bordism is $ 2 $-primary \cite{anderson1967structure}, and so any extra torsion term is itself divisible by $ p $.

By our main result, since all KO-theory (resp. K-theory) characteristic numbers are indices of twisted Dirac (resp. spin$ ^c $ Dirac) operators, their values mod $ p $ are fully determined up to an integer constant dependent only on $ M $, and hence so is the bordism class mod $ p $ by above. It therefore suffices to find an example in each case of such an $ M $ satisfying the conditions of the corollary and which isn't already divisible by $ p $ in bordism; the attached $ M $-dependent constants can then be immediately determined by checking any non-zero genus of each example.

For $ m<p-1 $ in both spin and complex cases, the result is immediate anyway by our initial discussion. For $ m=p-1 $ in the complex case, take $ \mathbb{C}P^{p-1} $, with a $ C_p^2 $ action where one generator acts on the projectivized coordinates cyclically and the other acts as diagonal multiplication with eigenvalues $ 1, \zeta, \zeta^2, \cdots, \zeta^{p-1} $. For $ m=p-1 $ in the spin case, take $ V_p $ represented by the Fermat hypersurface $ x_0^{2p}+\cdots+x_p^{2p}=0 $ in $ \mathbb{C}P^p $ and let $ (a, b)\in C_p^2 $ act on $ \mathbb{C}P^p $ as
\[ (a, b)\cdot [x_0, \cdots, x_p]=[x_0, \zeta^ax_1, \cdots, \zeta^{(p-1)a}x_{p-1}, \zeta^bx_p]; \]
this action then restricts to an action on $ V_p $ which has no stationary points. The appearance of $ 2^{-1} $ in the final result comes from $ \hat{A}(V_p)=2 $ (cf. $ \text{Td}(\mathbb{C}P^{p-1})=1 $).

Finally, we comment post facto on why we don't need an extra condition of the action of $ P $ lifting to the spin structure. Since $ M $ is connected and $ P $ has odd order, every element of $ P $ preserves the orientation of $ M $. Hence $ P $ acts on the set of possible spin structures on $ M $ over that orientation. This set forms a torsor over $ H^1(M;\mathbb{Z}/2\mathbb{Z}) $, and hence has size a power of $ 2 $. If it's just $ 1 $, it's certainly fixed by $ P $. If it's a larger power of $ 2 $, their $ P $-orbit sizes come in powers of $ p $; hence at least one $ P $-orbit has size $ 1 $, and there's some $ P $-invariant spin structure. A change of the spin structure changes the spin bordism class by a torsion element, which in spin bordism is $ 2 $-primary \cite{anderson1967structure} and therefore invisible mod $ p $, so we can assume the spin structure is this $ P $-invariant one. Moreover, the obstruction to equipping this invariant spin structure with a lift of the $ P $-action to the spinor bundle lives in $ H^2(P;\mathbb{Z}/2\mathbb{Z}) $, which vanishes since $ P $ is of odd order; see \cite{atiyah1970spin}, see also \cite[p.168]{farsi1992equivariant}. Therefore, we can take the action of $ P $ to lift to the spinor bundle, where it necessarily commutes with Clifford multiplication with respect to a $ P $-invariant metric and hence the Dirac operator.
\end{proof}

\newpage

\begin{corollary}\label{thm:main-thm-for-spin}
Let $ P $ be a $ p $-group for odd $ p $ which acts smoothly and without stationary points on a closed connected spin manifold $ M $ of dimension $ 2m $.
\begin{enumerate}[label=(\alph*)]
\item If $ m<p-1 $, then $ \hat{A}(M;W)\equiv 0\inmod{p} $ for every $ P $-equivariant complex vector bundle $ W $.
\item If $ m=p-1 $, then $ \hat{A}(M;W)\equiv\mathrm{rk}(W)\hat{A}(M)\inmod{p} $ for every $ P $-equivariant complex vector bundle $ W $ on which the kernel of $ P\to\mathrm{Diff}(M) $ acts trivially. In particular, $ \mathrm{sign}(M)\equiv\hat{A}(M)\inmod{p} $.
\end{enumerate}
\end{corollary}

\begin{proof}
By an identical line of reasoning as the discussion at the end of the proof of Corollary \ref{thm:divisibility-of-characteristic-numbers}, we can take the equivariance hypothesis of our main theorem to be satisfied for the spin Dirac operator without any further assumptions.

This corollary is then a direct application of the main theorem. To show $ \text{sign}(M)\equiv\hat{A}(M)\inmod{p} $ when $ m=p-1 $, recall the signature operator is $ D\otimes\mathcal{S}^* $ for $ D $ the Dirac operator and $ \mathcal{S} $ the spinor bundle, which has rank $ 2^{p-1} $. The result then follows by Fermat, $ 2^{p-1}\equiv 1\inmod{p} $.
\end{proof}

The following lemma allows us to modify the condition on the action of $ P $ in the main theorem from `has no stationary points' to `is effective' when $ P $ is non-abelian.

\begin{lemma}\label{thm:corollary-lemma}
Let $ P $ be a non-abelian $ p $-group for odd $ p $ which acts smoothly and effectively on a connected manifold $ M $ of dimension $ n $. If $ n<2p $, then the action of $ P $ has no stationary points.
\end{lemma}
\begin{proof}
We show the contrapositive: if the action has a stationary point, then $ n\ge 2p $. Suppose $ x\in M $ is such a stationary point; then $ V: P\to GL(T_xM) $ gives a real dimension $ n $ representation. We may assume that $ P $ acts on $ M $ by isometries; since isometries of a connected Riemannian manifold are determined by their value and derivative at one point \cite[Prop. 5.22]{lee2018introduction}, effectiveness of the $ P $ action implies that $ V $ is faithful. It then suffices to show that any real faithful representation of $ P $ has dimension at least $ 2p $.

Let $ V_\mathbb{C} $ be the complexification, which is still faithful and has complex dimension $ n $. A faithful action by a sum of one-dimensional representations would imply $ P $ were abelian, which we assume it isn't, so $ V_\mathbb{C} $ has an irreducible summand $ U $ with dimension greater than $ 1 $, and hence at least $ p $, as degrees of irreducible summands divide $ |P| $.

By Burnside, there are no non-trivial irreducible real-valued characters, so the conjugate representation $ \overline{U} $ is different from $ U $, and must also be contained in $ V_\mathbb{C} $ since the latter is conjugation-invariant. The summands $ U $ and $ \overline{U} $ together thus contribute dimension at least $ p+p=2p $, so $ n\ge 2p $.
\end{proof}

\begin{corollary}\label{thm:moonshine}
No closed, connected $ 24 $-dimensional spin manifold $ M $ with $ \hat{A}(M)=1 $ and $ \hat{A}(M;T_\mathbb{C}M)=0 $ can be acted on smoothly and non-trivially by the monster group.
\end{corollary}

\begin{proof}
Suppose such an action on such a manifold exists. The set of group elements of the monster acting trivially on $ M $ is a normal subgroup; therefore, since the monster is simple, any non-trivial action must be effective. Let $ p=13 $ and $ P $ be the Sylow $ 13 $-subgroup of the monster, which therefore also acts effectively on $ M $; this is an extraspecial group of order $ 13^{1+2} $ \cite{CCNPW85}.

By Lemma \ref{thm:corollary-lemma}, since $ P $ is non-abelian, the action of $ P $ on $ M $ has no stationary points. Part (b) of Corollary \ref{thm:main-thm-for-spin} therefore applies; taking $ W=T_\mathbb{C}M $, we get that $ \hat{A}(M;T_\mathbb{C}M)\equiv 24\hat{A}(M)\inmod{13} $, not true for the desired values of these quantities.
\end{proof}

This resolves a question of Hirzebruch. The author can't help but point out the traditional association of the number $ 13 $ with bad luck.

\bibliographystyle{amsplain}
\bibliography{references}

\end{document}